\documentclass[11pt]{article}

\usepackage{latexsym}
\usepackage{amsfonts}
\usepackage{enumerate}
\usepackage{multicol}
\usepackage{graphicx}
\usepackage{amssymb}
\usepackage{amsmath}
\usepackage{epic}

\usepackage{amsthm}

\begin{document}
\parindent=0.2in
\parskip 0in

\begin{flushright}

{\huge {\bf Report on the 63rd Annual International Mathematical Olympiad}}

\vspace*{.2in}

{\Large B\'ELA BAJNOK} \\

{\small Gettysburg College \\  Gettysburg, PA 17325 \\ {\tt bbajnok@gettysburg.edu}}

\end{flushright}

\vspace*{.2in}

\noindent The International Mathematical Olympiad (IMO) is the world's leading mathematics competition for high school students and is organized annually by different host countries.  The competition consists of three problems each on two consecutive days, with an allowed time of four and a half hours both days.  In recent years, more than one hundred countries have sent teams of up to six students to compete.  

The 63rd IMO was organized by Norway, and it was held in Oslo between July 6 and July 16, 2022, with the participation of 589 contestants from 104 countries.    

Each year, the members of the US team are chosen during the Math Olympiad Program (MOP), a year-long endeavor organized by the MAA's American Mathematics Competitions (AMC) program.  Students gain admittance to MOP based on their performance on a series of examinations, culminating in the USA Mathematical Olympiad (USAMO).  A report on the 2022 USAMO can be found in the February 2023 issue of this {\em Magazine}; a similar report on the 2022 USA Junior Mathematical Olympiad appeared in the January 2023 issue of the {\em College Mathematics Journal}.  More information on the American Mathematics Competitions program can be found on the site {\tt https://www.maa.org/math-competitions}.  

The members of the 2022 US team were Kevin Cong (12th grade, Phillips Exeter Academy, NH);  Ram Goel (12th grade, Krishna Homeschool, OR); Andrew Gu (12th grade, Tesoro High School, CA); Derek Liu (11th grade, Torrey Pines High School, CA); Luke Robitaille (12th grade, Robitaille Homeschool, TX); and  Eric Shen (11th grade, Lynbrook High School, CA).  Cong, Liu, Robitaille, and Shen each earned Gold Medals; Gu received a Silver Medal; and Goel earned a Bronze Medal.  In the unofficial ranking of countries, the United States finished third after China (first) and South Korea (second).

Below we present the problems and solutions of the 63rd IMO.  Our solutions are those of the current author, utilizing some of the various sources already available.  Each problem was worth 7 points; the nine-tuple $(a_7, a_6, a_5, a_4, a_3, a_2, a_1, a_0; {\bf a})$ states the number of students who scored $7, 6, \ldots , 0$ points, respectively, followed by the mean score achieved for the problem.  The combined mean total score of all participants was 20.33, substantially more than last year's 11.59.

\bigskip

\noindent {\bf Problem 1} $(385, 51, 13, 10, 21, 19, 56, 34; {\bf 5.54})$; {\em proposed by France.}  The Bank of Oslo issues two types of coin: aluminum (denoted A) and bronze (denoted B). Marianne has $n$ aluminum coins and $n$ bronze coins, arranged in a row in some
arbitrary initial order. A {\em chain} is any subsequence of consecutive coins of the same type. Given a
fixed positive integer $k \le 2n$, Marianne repeatedly performs the following operation: she identifies
the longest chain containing the $k$-th coin from the left, and moves all coins in that chain to the left end
of the row. For example, if $n = 4$ and $k = 4$, the process starting from the ordering $AABBBABA$
would be
$$AAB\underline{B}BABA \rightarrow BBB\underline{A}AABA \rightarrow AAA\underline{B}BBBA \rightarrow BBB\underline{B}AAAA \rightarrow BBB\underline{B}AAAA \rightarrow \cdots.$$
Find all pairs $(n, k)$ with $1 \leq k \leq 2n$ such that for every initial ordering, at some moment during
the process, the leftmost $n$ coins will all be of the same type.

\bigskip

\noindent  {\em Solution.}  Our claim is that the pairs satisfying the required property are those where $n \leq k \leq \lceil 3n/2 \rceil$.

We write our sequence of coins using exponential notation; for instance, $AABBBABA$ used in the example above will be written as $A^2B^3ABA$.  We also call a maximal chain of coins of the same type a {\em block}; for example, $A^2B^3ABA$ consists of five blocks.  The desired ordering of the $2n$ coins is then one of the two orderings that consist of only two blocks, $A^nB^n$ or $B^nA^n$.

To see that we must have $k \geq n$, observe that (for example) $A^{n-1}B^nA$ is unchanged by the operation when $k<n$.  We can rule out $k>\lceil 3n/2 \rceil$ as well: performing the operation repeatedly on $A^{\lceil n/2 \rceil} B ^{\lceil n/2 \rceil}A^{\lfloor n/2 \rfloor} B ^{\lfloor n/2 \rfloor}$ will leave the number of blocks at four since each step will simply move the last block to the front. 

Suppose now that $n \leq k \leq \lceil 3n/2 \rceil$, and consider an arbitrary ordering of the $2n$ coins.  If we have only two blocks, then we are done, so let's assume that our ordering contains at least three blocks.  Note that Marianne's operation never increases the number of blocks; our goal is to show that, performing the operation three times, the number of blocks will decrease by at least one.  This then means that the procedure eventually reaches the desired outcome: consisting of only two blocks.     

Observe that, unless the $k$-th coin is in the left-most block or the right-most block, moving the block containing it reduces the number of blocks by at least one, since the two neighboring blocks are then merged.  Now the $k$-th coin cannot be in a leftmost block, since $k \geq n$ and our sequence has at least three blocks.  Let the length of the rightmost block be $k_1$, the block before that have length $k_2$, and so on.  If performing the first operation does not reduce the number of blocks, then the $k$-th coin must be in a rightmost block, so $k \geq 2n-k_1+1$.  Furthermore, we have at least four blocks, since otherwise, moving the rightmost block to the front would merge it with what was originally the third block from the right as these two blocks contain coins of the same type.  So if the second and third operations don't reduce the number of blocks either, then we would have $k \geq 2n-k_2+1$ and $k \geq 2n-k_3+1$.  But this cannot be, since then we would have 
$$2 \lceil 3n/2 \rceil \geq 2k \geq (2n-k_1+1)+(2n-k_3+1),$$ implying $k_1+k_3 > n$, which is impossible as the corresponding blocks consist of coins of the same type.

\bigskip

\noindent {\bf Problem 2} $(303, 19, 3, 7, 41, 23, 89, 104; {\bf 4.31})$; {\em proposed by the Netherlands.}  Let $\mathbb{R}^+$
 denote the set of positive real numbers. Find all functions $f : \mathbb{R}^+ \rightarrow \mathbb{R}^+$ such
that for each $x \in \mathbb{R}^+$, there is exactly one $y \in \mathbb{R}^+$ satisfying
$$xf(y)+yf(x) \leq 2.$$

\bigskip

\noindent  {\em Solution.}  We prove that the only function satisfying this property is $f(x)=1/x$.  

First, we observe that this function indeed works, since by the AM--GM inequality, for each $x \in \mathbb{R}^+$, the unique $ y \in \mathbb{R}^+$ for which $x/y + y/x \leq 2$ holds is $y=x$.

Suppose now that $f : \mathbb{R}^+ \rightarrow \mathbb{R}^+$ is a function satisfying the requirement that each $x \in \mathbb{R}^+$ has a unique {\em partner} $ y \in \mathbb{R}^+$ such that $$xf(y)+yf(x) \leq 2.$$  Note that, by symmetry, this also means that $x$ is the partner of $y$.

We claim that $f$ must be strictly decreasing.  Indeed, suppose that $x_1, x_2 \in \mathbb{R}^+$ with $x_1 < x_2$, and let $y_2$ be the partner of $x_2$.  Since $x_1$ is then not the partner of $y_2$ but $x_2$ is, we must have $$x_2f(y_2)+y_2f(x_2) \leq 2 < x_1f(y_2)+y_2f(x_1).$$  Since $x_1 < x_2$, we get $f(x_2) < f(x_1)$, proving our claim.

We can now establish that the partner of each positive real number is itself.  For that, suppose that $x \in \mathbb{R}^+$ has partner $y \in \mathbb{R}^+$ and that $x \neq y$.  Then $x$ is not its own partner, so $x f(x) + x f(x)$ is greater than $2$, and so $xf(x)>1$; similarly, $yf(y)>1$.  Adding these two inequalities to $(x-y)(f(y)-f(x))$, which is positive since $f$ is strictly decreasing, we get
$$ xf(y)+yf(x) = (x-y)(f(y)-f(x)) +xf(x)+yf(y)  > 2,$$ contradicting the fact that $x$ and $y$ are partners.       

Finally, we show that $f(x)=1/x$ for every $x \in \mathbb{R}^+$.  Since every $x \in \mathbb{R}^+$ is its own partner, we have $xf(x) \leq 1$.  Suppose, indirectly, that $f(x) < 1/x$ holds for some $x \in \mathbb{R}^+$, and let $z=1/f(x)$. Since $z > x$ and $f$ is strictly decreasing, we have $f(z) < f(x)$, and thus
$$z f(x) + x f(z) < 1 + x f(x) \leq 2,$$ but that would mean that $x$ and $z$ are partners, which is a contradiction.  Therefore, $f(x)=1/x$ for every $x \in \mathbb{R}^+$, as claimed.   

\bigskip

\noindent {\bf Problem 3} $(28, 4, 4, 3, 6, 69, 68, 407; {\bf 0.81})$; {\em proposed by the United States of America.}  Let $k$ be a positive integer and let $S$ be a finite set of odd prime numbers. Prove that there is at most one way (up to rotation and reflection) to place the elements of $S$ around a circle so that the product of any two neighbors is of the form $x^2 + x + k$ for some positive integer $x$.

\bigskip

\noindent  {\em Solution.}  Let $k$ be a fixed positive integer.  We will say that two odd primes are {\em $k$-friends} if their product is of the form $x^2 + x + k$ for some {\em nonnegative} integer $x$.  We will prove the slightly more general statement that, given a finite set $S$ consisting of odd primes, there is at most one way (up to rotation and reflection) to place the elements of $S$ around a circle so that any two neighbors are $k$-friends.  We will say that such an arrangement of the elements of $S$ is {\em valid}.  

We use induction on the size of $S$, noting that the claim is trivially true if $|S| \leq 3$.  In particular, what we will show is that if $p$ is the greatest prime in $S$, then from each valid arrangement of the elements of $S$ we get a valid arrangement of the elements of $S \setminus \{p\}$ by simply removing $p$.  It is easy to see that to establish this, it suffices to prove that $p$ may have at most two $k$-friends, and if $p$ has two $k$-friends, then these elements are also $k$-friends with each other.

We first rule out the possibility of $p$ having more than two $k$-friends that are less than $p$.  Suppose that $p$ has $k$-friends $q$, $r$, and $s$ with $p>q>r>s$, and let
$$pq=x^2 + x + k, \; \; pr=y^2+y+k, \; \; ps=z^2+z+k$$
for some nonnegative integers $x$, $y$, and $z$.  Since $p$ is larger than $q$, $r$, and $s$, we see that $x$, $y$, and $z$ are each less than $p$.  Taking the difference of the first two equations, we get 
$$p(q-r)=(x-y)(x+y+1),$$ and thus $x-y$ or $x+y+1$ is divisible by $p$.  But $1 \leq x- y \leq p-1$ and $1 \leq x+y+1 \leq 2p-1$, so in fact $x+y+1=p$.  Similarly, $x+z+1=p$ and $y+z+1=p$, so $x=y=z$, a contradiction.  

Suppose now that $p$ has $k$-friends $q$ and $r$.  Using the notations from above, we see that $p=x+y+1$ and thus $q-r=x-y$.    
We now verify that $q$ and $r$ are $k$-friends, as follows:
\begin{eqnarray*}
qr & = & pq -pq + qr \\
& = & (x^2+x+k) - (x+y+1)q + q (q-x+y) \\
& = & (x-q)^2 +(x-q) + k,
\end{eqnarray*}
and this completes our proof when $x-q \geq 0$.  But when $x-q < 0$, then $q-x-1 \geq 0$, in which case we can write
$$qr = (x-q)^2 +(x-q) + k = (q-x-1)^2 +(q-x-1) + k.$$ Now our proof is complete.

\bigskip

\noindent {\bf Problem 4} $(433, 1, 5, 10, 17, 18, 31, 74; {\bf 5.47})$; {\em proposed by Slovakia.}  Let $ABCDE$ be a convex pentagon such that $BC = DE$. Assume that there is a point $T$ inside $ABCDE$ with $TB = TD$, $TC = TE$, and $\angle ABT = \angle TEA$. Let line $AB$ intersect lines $CD$ and $CT$ at points $P$ and $Q$, respectively. Assume that the points $P, B, A, Q$ occur on their line in that order. Let line $AE$ intersect lines $CD$ and $DT$ at points $R$ and $S$, respectively. Assume that the points $R, E, A, S$ occur on their line in that order. Prove that the points $P, S, Q, R$ lie on a circle.

\bigskip

\noindent  {\em Solution.}  We start by observing that, according to our assumptions, triangles $BCT$ and $DET$ are congruent, and thus $\angle CTB = \angle ETD$.  (An illustration for this problem is provided in Figure 1.)  

\begin{figure}[ht]
\centering
\includegraphics{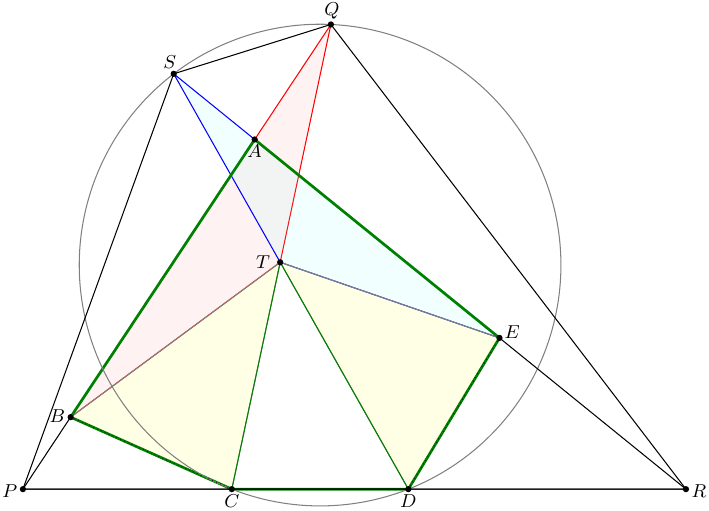}
\caption{Illustration for Problem 4}
\end{figure}

Then for the corresponding supplementary angles we have $\angle BTQ = \angle STE$; since $\angle ABT = \angle TEA$ as well, triangles $QBT$ and $SET$ are similar, and thus $\angle TQB = \angle EST$.   Furthermore, we have $$  \frac{TE}{TB} = \frac{TS}{TQ},$$ from which
$$TQ \cdot TE = TS \cdot TB.$$ But $TE = TC$ and $TB=TD$, so $$TQ \cdot TC = TS \cdot TD.$$  The converse of the Power of a Point Theorem then implies that $CDQS$ is cyclic.  Therefore, $\angle CQS = \angle CDS$, and thus 
$$\angle BQS = \angle CQS - \angle TQB = \angle CDS - \angle EST = 180^\circ - \angle SDR - \angle RSD =\angle DRS.   $$  This means that $\angle PQS = \angle PRS$, and thus $PRQS$ is cyclic, as claimed.

\bigskip

\noindent {\bf Problem 5} $(171, 49, 35, 15, 48, 44, 115, 112; {\bf 3.52})$; {\em proposed by Belgium.}  Find all triples $(a,b,p)$ of positive integers with $p$ prime and\[a^p = b! + p.\]

\bigskip

\noindent  {\em Solution.}  We claim that there are two such triples: $(a,b,p)=(2,2,2)$ and $(a,b,p)=(3,4,3)$.  It is clear that our equation is satisfied with these two choices; we will show that there are no others.

Suppose, then, that positive integers $a$, $b$, and $p$ satisfy the equation and that $p$ is prime.   
We first prove that, unless $a=b=p=2$, $b$ is strictly between $p$ and $2p$.  Indeed, if we were to have $b \geq 2p$, then $b!+p$ would be divisible by $p$ but not by $p^2$, and this is never the case for $a^p$.  On the other hand, if $b \leq p$, then we would also have $a \leq b$, since otherwise
$$a^p \geq (b+1)^p \geq b^p + pb+1 > b^b +p \geq b!+p,$$ a contradiction.  So $1 < a \leq b$, and thus $b!$ is divisible by $a$, but then so is $p=a^p-b!$, which can only happen if $a=p$.  This then means that 
$p^p=b!+p$, so by our assumption, $p^p \leq p!+p$, but that is possible  only for $p=2$ (and thus $a=b=p=2$), since otherwise
$$p! + p \leq p^{p-1} + p \leq 2p^{p-1} < p^p.$$

Next, we prove that $a=p$.  Since $b \geq p$ by the previous paragraph, $b!$, and therefore $a^p=b!+p$, must be divisible by $p$, and thus $a=qp$ for some positive integer $q$.  We need to show that $q=1$.  If $q < p$, then also $q < b$, so $b!$ and therefore $p=a^p-b!$ is divisible by $q$, and thus indeed $q=1$. Let us consider now the possibility that $q \geq p$, in which case $a \geq p^2$, and thus $b!+p = a^p \geq p^{2p}$.  On the other hand, since $b \leq 2p-1$, we have $b!+p \leq (2p-1)!+p$.  Therefore, we get that $(2p-1)!+p$ is at least $p^{2p}$, which is not possible as we now show.  Let us write
$$(2p-1)! = 1 \cdot (2p-1) \cdot 2 \cdot (2p-2) \cdot \cdots \cdot (p-1) \cdot (p+1) \cdot p.$$  Observe that for each $1 \leq i \leq p-1$, $i \cdot (2p-i)$ is less than $p^2,$  so $(2p-1)!$ is less than $p^{2p-1}$, and thus $(2p-1)!+p$ is less than $p^{2p}$, as claimed.

We can now establish our claim by proving that the equation $p^p=b!+p$ only holds for $(p,b)=(2,2)$ and $(p,b)=(3,4)$.  For the sake of a contradiction, assume that a solution $(p,b)$ exists for some prime $p \geq 5$.  Note that the integers $2, (p+1)/2$, and $p+1$ are pairwise distinct in this case, and since $b > p$, their product must be a divisor of $b!$.  We will arrive at a contradiction by proving that $p^p-p$ is never divisible by $(p+1)^2$.  
To see this, we can write
$$p^p-p=p \cdot(p^2-1)\cdot \left  (p^{p-3}+p^{p-5}+\cdots + p^2+1 \right),$$ so our claim follows if 
$$p \cdot(p-1)\cdot \left  (p^{p-3}+p^{p-5}+\cdots + p^2+1 \right)$$ is not divisible by $p+1$.  And indeed it isn't, since it is congruent to 
$$(-1) \cdot (-2) \cdot (p-1)/2 = p-1$$ modulo $p+1$.  This completes our proof.

\bigskip

\noindent {\bf Problem 6} $(22, 2, 3, 19, 7, 22, 80, 434; {\bf 0.68})$; {\em proposed by Serbia.}  Let $n$ be a positive integer. A {\em Nordic square} is an $n \times n$ board containing all the integers from $1$ to $n^2$ so that each cell contains exactly one number. Two different cells are considered {\em adjacent} if they share an edge. Every cell that is adjacent only to cells containing larger numbers is called a {\em valley}. An {\em uphill path} is a sequence of one or more cells such that
\begin{itemize}
  \item  the first cell in the sequence is a valley,
  \item each subsequent cell in the sequence is adjacent to the previous cell, and
  \item the numbers written in the cells in the sequence are in increasing order.
\end{itemize}
Find, as a function of $n$, the smallest possible total number of uphill paths in a Nordic square.

\bigskip

\noindent  {\em Solution.}  We claim that the smallest possible number of uphill paths is $2n^2-2n+1$.

To see that this is a lower bound, first observe that every Nordic square has at least one valley (for example, the cell containing $1$) and thus has at least one uphill path consisting of a single cell.  Now consider two adjacent cells, and suppose that the numbers in these cells are $k_1$ and $k_2$ with $k_1 > k_2$.  Note that one can always extend the path $k_1, k_2$ to a {\em downhill path} $k_1, k_2, \ldots, k_t$ so that each number in the sequence (except $k_1$) is less than the one before it, and that $k_t$ is a valley (here $2 \leq t \leq n^2$).  Therefore, the reverse path is an uphill path ending in the two adjacent cells we chose.  Since the board has $2n(n-1)$ pairs of adjacent cells, we have at least this many uphill paths containing two or more cells.  This proves that the number of uphill paths is at least $2n(n-1)+1$.      

Our remaining task is to provide an example of an $n \times n$ Nordic square (for every $n$) that has exactly one valley and where no two adjacent cells are the last two cells of two different uphill paths.  We will accomplish this via the following lemma.

\noindent {\bf Lemma.}  {\em For any pair of positive integers $(m,n)$ it is possible to mark some of the cells of an $m \times n$ board so that no two marked cells are adjacent and that the unmarked cells form a tree (that is, for any pair of unmarked cells there is a unique path  connecting them that avoids all marked cells).}

\noindent {\em Proof of Lemma.}  For $n=1$ the claim is obvious, so let us assume that $n \geq 2$.  We first provide an example for the case $m=6$, as follows: we mark cells $(1,j)$ and $(3,j)$ for all even $j$; mark cells $(4,j)$ and $(6,j)$ for all odd $j >1$; and mark the cell $(5,1)$.  (An illustration of this for $n=13$ can be seen in Figure 2.)  Then, to get a construction for $m=6q$ for any positive integer $q$, we juxtapose $q$ copies of the $6 \times n$ board, and from this we get the case of $m=6q-r$ with $1 \leq r \leq 5$ by deleting the first $r$ rows when $r \in \{1,3,4\}$, while deleting the first $r-1$ rows and the last row when $r \in \{2,5\}$.  It is not hard to verify that this construction works.

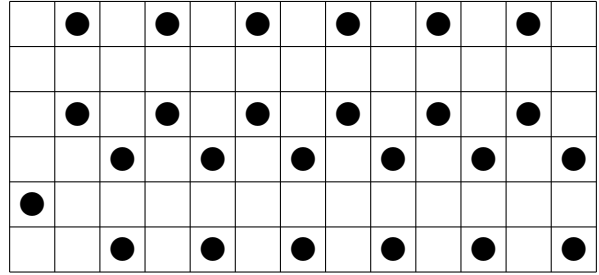
\begin{figure}[h]
\begin{center}
\setlength{\unitlength}{.6cm}
\begin{picture}(13,6)
\put(0,0){\line(1,0){13}}
\put(0,1){\line(1,0){13}}
\put(0,2){\line(1,0){13}}
\put(0,3){\line(1,0){13}}
\put(0,4){\line(1,0){13}}
\put(0,5){\line(1,0){13}}
\put(0,6){\line(1,0){13}}

\put(0,0){\line(0,1){6}}
\put(1,0){\line(0,1){6}}
\put(2,0){\line(0,1){6}}
\put(3,0){\line(0,1){6}}
\put(4,0){\line(0,1){6}}
\put(5,0){\line(0,1){6}}
\put(6,0){\line(0,1){6}}
\put(7,0){\line(0,1){6}}
\put(8,0){\line(0,1){6}}
\put(9,0){\line(0,1){6}}
\put(10,0){\line(0,1){6}}
\put(11,0){\line(0,1){6}}
\put(12,0){\line(0,1){6}}
\put(13,0){\line(0,1){6}}

\put(1.5,3.5){\circle*{0.5}}
\put(3.5,3.5){\circle*{0.5}}
\put(5.5,3.5){\circle*{0.5}}
\put(7.5,3.5){\circle*{0.5}}
\put(9.5,3.5){\circle*{0.5}}
\put(11.5,3.5){\circle*{0.5}}

\put(1.5,5.5){\circle*{0.5}}
\put(3.5,5.5){\circle*{0.5}}
\put(5.5,5.5){\circle*{0.5}}
\put(7.5,5.5){\circle*{0.5}}
\put(9.5,5.5){\circle*{0.5}}
\put(11.5,5.5){\circle*{0.5}}

\put(2.5,0.5){\circle*{0.5}}
\put(4.5,0.5){\circle*{0.5}}
\put(6.5,0.5){\circle*{0.5}}
\put(8.5,0.5){\circle*{0.5}}
\put(10.5,0.5){\circle*{0.5}}
\put(12.5,0.5){\circle*{0.5}}

\put(2.5,2.5){\circle*{0.5}}
\put(4.5,2.5){\circle*{0.5}}
\put(6.5,2.5){\circle*{0.5}}
\put(8.5,2.5){\circle*{0.5}}
\put(10.5,2.5){\circle*{0.5}}
\put(12.5,2.5){\circle*{0.5}}

\put(0.5,1.5){\circle*{0.5}}

\end{picture}
\end{center} 
\caption{The markings on a $6 \times 13$ board}
\end{figure}

We now create our $n \times n$ Nordic square as follows.  Having the set of markings on the board from our Lemma, we place $1$ in an arbitrary unmarked cell, then fill in the rest of the unmarked cells one at a time by placing the smallest yet-unplaced integer adjacent to a cell that already has a number in it.  Since our tree is connected, this procedure will only terminate once all unmarked cells contain a number.  The rest of the numbers can then be placed in the marked cells arbitrarily.

We need to verify that this placement of the numbers yields exactly $2n^2-2n+1$ uphill paths.  In order to do so, it will suffice to prove that we have exactly one valley and that any two adjacent cells are the last two cells in a unique uphill path.  To see that we have exactly one valley (namely, the cell containing $1$), note that no other unmarked cell can be a valley since the number that it contains is larger than the number in its neighbor that was numbered earlier, and clearly none of the marked cells is a valley since all their neighbors contain smaller numbers as no two marked cells are adjacent.  The fact that no two adjacent cells may be the last two cells in two different uphill paths follows easily as well: the cell with the smaller number in it must be unmarked, and thus there is a unique path from it to the cell containing $1$ because the unmarked cells form a tree.  This completes our proof.

\bigskip
\noindent {\small {\bf Acknowledgments.}  I am very grateful to Evan Chen, Peter Francis, Jerrold Grossman, Jonathan Kane, and Enrique Trevi\~no for their willingness to proofread parts of this article and for doing it so carefully.  Credit for Figure 1 goes to Evan Chen.

\bigskip
\noindent {\small {\bf Summary.}  We present the problems and solutions to the 63rd Annual International Mathematical Olympiad.

\bigskip

\noindent {\bf B\'ELA BAJNOK} (MR Author ID: 314851) is a Professor of Mathematics at Gettysburg College and the Director of the American Mathematics Competitions program of the MAA.

\end{document}